\documentclass[11pt]{article}

\usepackage{latexsym}
\usepackage{amssymb}

\makeatletter \@addtoreset{equation}{section}

\newtheorem{theorem}{Theorem}[section]
\newtheorem{proposition}[theorem]{Proposition}
\newtheorem{corollary}[theorem]{Corollary}
\newtheorem{lemma}[theorem]{Lemma}

\newtheorem{remarks}[theorem]{Remarks}

\def\endverif{\nopagebreak\newline\mbox{\ }\hfill\rule{2mm}{2mm}}
\def\endveriff{\nopagebreak\mbox{\ }\hfill\rule{2mm}{2mm}}

\def\1{{\bf 1}}

\def\2{{1/2}}

\def\C{{\mathbb C}}
\def\N{{\mathbb N}}

\def\Z{{\mathbb Z}}

\def\A{{\cal A}}

\def\eps{\varepsilon}

\def\CPA{{\rm CPA}}

\def\Mat{{\rm Mat}}
\def\rank{{\rm rank}\,}

\begin{document}

\title{\bf The McMillan theorem for a class of asymptotically
abelian C$^*$-algebras}

\author{
Sergey Neshveyev$^{1)}$ \\
{\it\small B. Verkin Institute for Low
Temperature Physics and Engineering,}\\
%{\it\small National Academy of Sciences of Ukraine,}\\
{\it\small 47 Lenin Ave., Kharkov 310164, Ukraine}\\
\\
and\\
\\
Erling St{\o}rmer$^{2)}$ \\
{\it\small Department of Mathematics, University of Oslo,}\\
{\it\small P.O. Box 1053, Blindern, 0316 Oslo, Norway}}
\date{}

\footnotetext[1]{The Royal Society/NATO postdoctoral fellow.
Partially supported by Award No UM1-2092 of the Civilian Research
\& Development Foundation and NATO grant SA (PST.CLG.976206)5273.}

\footnotetext[2]{Partially supported by the Norwegian Research
Council.}

\maketitle

\begin{abstract}
An extension of the Shannon-McMillan-Breiman theorem to a class of
non-commutative dynamical systems is given.
\end{abstract}

%\bigskip

\section{Introduction}
In ergodic theory one of the main theorems on entropy is the
McMillan theorem, also called the Shannon-McMillan-Breiman
theorem. In one form it states that if $(X,\mu)$ is a probability
space and $T$ a measure preserving ergodic transformation, then
for any finite measurable partition $\xi$ and any $\eps>0$ there
exists $n_0$ such that if $n\ge n_0$ then outside of a union of
atoms of total measure $<\eps$ every atom in
$\vee^{n-1}_{i=0}T^{-i}\xi$ has measure in the interval
$(e^{-n(H(T;\xi)+\eps)},e^{-n(H(T;\xi)-\eps)})$.

In the present paper we shall give a non-commutative extension of
the McMillan theorem. Our setting will be asymptotically abelian
C$^*$-dynamical systems with locality $(A,\tau,\alpha)$~\cite{NS},
where we assume $\tau$ is an invariant ergodic trace. If $N$ is a
local subalgebra of $A$ which is a mean generator~\cite{GS2}, i.e.
the C$^*$-algebra $\vee_{i\in\Z}\alpha^i(N)=A$ and
$H(\alpha)=\lim{1\over n}H(\vee^{n-1}_{i=0}\alpha^i(N))$, then we
shall show that for any $\eps>0$ and any $n$ sufficiently large
their exists a central projection $z_n$ in
$N_n=\vee^{n-1}_{i=0}\alpha^i(N)$ such that $\tau(z_n)<\eps$ and
$$
e^{-n(H(\alpha)+\eps)}<\tau(e)<e^{-n(H(\alpha)-\eps)}
$$
for any minimal projection $e$ in $N_n(1-z_n)$. The proof is based
on ideas in~\cite{NS} and uses the classical McMillan theorem.

The above result forms the main contents of Section~\ref{2}. Then
in Section~\ref{3} we study sufficient conditions for a local
algebra to be a mean generator. For example, we shall show that if
$N$ is a mean generator and $M$ is a C$^*$-algebra contained in
$\vee^n_{-n}\alpha^i(N)$ then $M$ is itself a mean generator for
the C$^*$-algebra $\vee_{i\in\Z}\alpha^i(M)$.

\bigskip\bigskip

\section{The McMillan Theorem} \label{2}

Throughout the paper we consider C$^*$-dynamical systems
$(A,\tau,\alpha)$ which are asymptotically abelian with
locality~\cite{NS}, i.e. we assume that there exists a dense
$\alpha$-invariant $*$-subalgebra $\A$ of~$A$ such that for all
pairs $a,b\in\A$ the C$^*$-algebra they generate is finite
dimensional, and there is $p=p(a,b)\in\N$ such that
$[\alpha^j(a),b]=0$ for $|j|\ge p$. Recall that by a local algebra
we mean a finite dimensional subalgebra of $\A$. We shall assume
also that $\tau$ is an $\alpha$-invariant trace.

\begin{theorem} \label{2.1}
Suppose the trace $\tau$ is ergodic, i.e. is an extremal
$\alpha$-invariant state. Suppose $N$ is a local algebra  which
is a mean generator and the entropy $H(\alpha)$ is finite, so
$$
H(\alpha)=\lim_{n\to\infty}{H(N_n)\over n}<\infty,
$$
where $N_n=\vee^{n-1}_{i=0}\alpha^i(N)$. Then given $\eps>0$ there
exist $n_0$ and central projections $z_n$ in $N_n$ such that
$\tau(z_n)<\eps$ for $n\ge n_0$ and
$$
e^{-n(H(\alpha)+\eps)}<\tau(e)<e^{-n(H(\alpha)-\eps)}
$$
for any minimal projection $e$ in $N_n(1-z_n)$.
\end{theorem}

Note that since $\tau$ is a trace, ergodicity of $\tau$ is
equivalent to ergodicity of the automorphism~$\alpha$ on the weak
closure $\pi_\tau(A)''$ of the GNS-representation of $A$. Note
also that we can consider $\pi_\tau(A)$ instead of $A$, and so
assume that $\tau$ is faithful.

\medskip

For all examples we know, the assumption on mean generator is
always fulfilled. However, we were unable to prove that it is
automatic. In the next section we shall discuss several sufficient
conditions for a local algebra to be a mean generator.

\medskip

To prove Theorem~\ref{2.1} we shall need the first part of the
following proposition.

\begin{proposition} \label{2.2}
Let $N$ be a local algebra, $p\in\N$ such that $\alpha^j(N)$
commutes with $N$ for $|j|\ge p$. Then

{\rm(i)} $\displaystyle
H(\alpha|_{A(N)})=\lim_{n\to\infty}{H(N_{n-p}\,,\alpha^n)\over
n},$ where $A(N)=\vee_{i\in\Z}\alpha^i(N);$

{\rm(ii)} $\displaystyle hcpa_\tau(\alpha|_{A(N)})
=\lim_{n\to\infty}{hcpa_\tau(N_{n-p}\,,\alpha^n)\over n},$ where
$hcpa_\tau(N_{n-p}\,,\alpha^n)$ is Voiculescu's completely
positive approximation entropy of $\alpha^n$ computed with respect
to any finite set spanning $N_{n-p}$.
\end{proposition}

\noindent{\it Proof.} (i) The proof is similar to
\cite[Lemma~3.4]{NS}. Since $H(N_{n-p}\,,\alpha^n)\le
nH(\alpha|_{A(N)})$, it is enough to prove that for any $m\in\N$
$$
H(N_m\,,\alpha)\le\liminf_{n\to\infty}{H(N_{n-p}\,,\alpha^n)\over
n}.
$$
For this fix $m_0$ and take $n\ge m+m_0+p$ and $k\in\N$. Consider
the algebra $B=(\vee^k_{j=0}\alpha^{jn}(N_{n-p}))^{m_0}$ and
define unital completely positive mappings $\psi\colon B\to A$ and
$\rho\colon A\to B$ as follows:
$$
\psi(b_1,\ldots,b_{m_0})={1\over m_0}
   \sum^{m_0}_{i=1}\alpha^{-i+1}(b_i),
$$
$$
\rho(a)=(E(a),(E\circ\alpha)(a),\ldots,(E\circ\alpha^{m_0-1})(a)),
$$
where $E\colon A\to\vee^k_{j=0}\alpha^{jn}(N_{n-p})$ is the
$\tau$-preserving conditional expectation. For $i=0,\ldots,kn-1$
and $a\in N_m$, we have from the proof of \cite[Lemma~3.4]{NS}
$$
||(\psi\circ\rho)(\alpha^i(a))-\alpha^i(a)||\le{2p\over m_0}||a||.
$$
Let $\theta\colon N_m\hookrightarrow A$ be the inclusion mapping.
Then by \cite[Theorem~VI.3 and Proposition~III.6]{CNT}
\begin{eqnarray*}
{1\over kn}
H_\tau(\theta,\alpha\circ\theta,\ldots,\alpha^{kn-1}\circ\theta)
&\le&{1\over kn}H_\tau(\psi\circ\rho\circ\theta,
    \psi\circ\rho\circ\alpha\circ\theta,\ldots,
    \psi\circ\rho\circ\alpha^{kn-1}\circ\theta)+\eps \\
 &\le&{1\over kn}H_\tau(\psi)+\eps \\
 &\le&{1\over kn}S(\tau\circ\psi)+\eps \\
 &=&{1\over kn}
     S(\tau|_{\vee^k_{j=0}\alpha^{jn}(N_{n-p})})+{\log m_0\over kn}+\eps,
\end{eqnarray*}
where $\eps=\eps(\rank N_m,p/m_0)\to0$ as $m_0\to\infty$. Since
$\alpha^{jn}(N_{n-p})$ commutes with $N_{n-p}$ for $j\ne0$,
$$
S(\tau|_{\vee^k_{j=0}\alpha^{jn}(N_{n-p})})
=H(N_{n-p}\,,\alpha^n(N_{n-p}),\ldots,\alpha^{kn}(N_{n-p})).
$$
Thus letting $k\to\infty$ in the above inequality we get
$$
H(N_m\,,\alpha)\le{1\over n}H(N_{n-p}\,,\alpha^n)+\eps,
$$
and the proof is complete.

(ii) This part is proved analogously: given
$(B_0,\psi_0,\rho_0)\in\CPA(A,\tau)$ which approximates
$\cup^k_{j=0}\alpha^{jn}(N_{n-p})$ we construct
$(B_1,\psi_1,\rho_1)$, where $B_1=B^{m_0}_0$,
$$
\psi_1(b_1,\ldots,b_{m_0})={1\over m_0}
   \sum^{m_0}_{i=1}\alpha^{-i+1}(\psi_0(b_i)),
$$
$$
\rho_1(a)=(\rho_0(a),(\rho_0\circ\alpha)(a),\ldots,
(\rho_0\circ\alpha^{m_0-1})(a)),
$$
which approximates $\cup^{kn-1}_{j=0}\alpha^j(N_m)$.
\endverif

Since $H_{\tau_1\otimes\tau_2}(M_1\otimes
M_2)=H_{\tau_1}(M_1)+H_{\tau_2}(M_2)$ for finite dimensional
C*-algebras $M_1$ and $M_2$ with normalized traces $\tau_1$ and
$\tau_2$ respectively, we have

\begin{corollary}
For systems which are asymptotically abelian with locality the
tensor product formula for the entropy with respect to tracial
states holds.
\end{corollary}

Let $\{M_n\}_n$ be a sequence of finite dimensional
C$^*$-algebras, $h$ a non-negative number, $\eps>0$. We shall say
that the McMillan theorem holds for $\{M_n\}_n$, $h$ and $\eps$,
if the conclusion of Theorem~\ref{2.1} holds for $M_n$ instead of
$N_n$, and $h$ instead of $H(\alpha)$. So there exist $n_0$ and
central projections $z_n$ in $M_n$ such that $\tau(z_n)<\eps$ for
$n\ge n_0$ and
$$
e^{-n(h+\eps)}<\tau(e)<e^{-n(h-\eps)}
$$
for any minimal projection $e$ in $M_n(1-z_n)$.

\begin{lemma} \label{2.4}
For any $\eps>0$ there exists $\eps_1>0$ such that if
$\{\xi_n\}_n$ and $\{\zeta_n\}_n$ are sequences of finite
partitions of a Lebesgue space $(X,\mu)$, $\zeta_n\prec\xi_n$,
$h=\lim{1\over n}H(\xi_n)$, $h_1=\lim{1\over n}H(\zeta_n)$,
$|h-h_1|<\eps_1$ and the McMillan theorem holds for
$\{\zeta_n\}_n$, $h_1$ and $\eps_1$, then it holds also for
$\{\xi_n\}_n$, $h$ and $\eps$.
\end{lemma}

\noindent{\it Proof.} Let $\{C_{ni}\}_{i\in X_n}$ be the atoms of
$\xi_n$, $\{D_{nj}\}_{j\in Y_n}$ the atoms of $\zeta_n$. Let
$\tilde Y_n$ be the set of $j\in Y_n$ for which
$$
e^{-n(h_1+\eps_1)}<\mu(D_{nj})<e^{-n(h_1-\eps_1)}.
$$
For $n$ large enough
$$
0\le{1\over
n}(H(\xi_n)-H(\zeta_n))<(h+\eps_1)-(h_1-\eps_1)<3\eps_1,
$$
or
\begin{equation} \label{e2.1}
0\le{1\over n}\sum_{i\in X_n\atop j\in Y_n}\mu(C_{ni}\cap
D_{nj})\log{\mu(D_{nj})\over\mu(C_{ni})}<3\eps_1.
\end{equation}
For $i\in X_n$ let $j(i)$ be a unique index such that
$C_{ni}\subset D_{nj(i)}$. Let
$$
X'_n=\{i\in X_n\ |\ j(i)\in\tilde Y_n\},\ \ X''_n=\{i\in X_n\ |\
\mu(D_{nj(i)})<e^{n\sqrt{\eps_1}}\mu(C_{ni})\},\ \ \tilde
X_n=X'_n\cap X''_n.
$$
By virtue of (\ref{e2.1})
$$
\sum_{i\in X_n\backslash X''_n}\mu(C_{ni})<3\sqrt\eps_1.
$$
Hence, for $n$ large enough,
$$
\sum_{i\notin\tilde X_n}\mu(C_{ni})\le\sum_{i\notin
X'_n}\mu(C_{ni})+\sum_{i\notin
X''_n}\mu(C_{ni})\le\sum_{j\notin\tilde
Y_n}\mu(D_{nj})+\sum_{i\notin
X''_n}\mu(C_{ni})<\eps_1+3\sqrt\eps_1.
$$
For $i\in\tilde X_n$
$$
\mu(C_{ni})\le\mu(D_{nj(i)})<e^{-n(h_1-\eps_1)}<e^{-n(h-2\eps_1)}
$$
and
$$
\mu(C_{ni})\ge
e^{-n\sqrt\eps_1}\mu(D_{nj(i)})>e^{-n(h_1+\eps_1+\sqrt\eps_1)}
>e^{-n(h+2\eps_1+\sqrt\eps_1)}.
$$
Thus we can take $\eps_1$ such that $2\eps_1+3\sqrt\eps_1<\eps$.
\endverif

\begin{lemma} \label{2.5}
Let $\{\xi_n\}_n$ be an increasing sequence of finite partitions
of a Lebesgue space $(X,\mu)$, $\zeta$ a finite partition such
that $\zeta\prec\vee_n\xi_n$. Suppose the McMillan theorem holds
for $\{\xi_n\vee\zeta\}_n$, $h$ and $\eps$. Then it holds also for
$\{\xi_n\}_n$, $h$ and $2\eps$.
\end{lemma}

\noindent{\it Proof.} Let $\{C_{ni}\}_{i\in X_n}$ be the atoms of
$\xi_n$, $\{D_j\}_{j\in Y}$ the atoms of $\zeta$. Let
$$
X'_n=\{i\in X_n\ |\ \exists\ j(i)\in Y\ :\ \mu(C_{ni}\cap
D_{j(i)})>(1-\eps/2)\mu(C_{ni})\}.
$$
We assert that
$$
\sum_{i\notin X'_n}\mu(C_{ni})\to0\ \ \hbox{as}\ \ n\to\infty.
$$
Indeed, otherwise there exist $j\in Y$ and $c>0$ such that
$$
\sum_{i\notin X'_n}\mu(C_{ni}\cap D_j)\ge c
$$
for infinitely many $n$'s. In other words, the set
$E_n=\cup_{i\notin X'_n}(C_{ni}\cap D_j)$ has measure $\mu(E_n)\ge
c$, and for the conditional expectation $E(D_j|\xi_n)$ of the
characteristic function of the set $D_j$ with respect to the
partition $\xi_n$ we have
$$
E(D_j|\xi_n)(x)\le1-{\eps\over2}\ \ \hbox{on}\ \ E_n,
$$
which contradicts the a.e. convergence of $\{E(D_j|\xi_n)\}_n$ to
the characteristic function of $D_j$.

Now let $Z_n$ be the set of pairs $(i,j)$ such that
$$
e^{-n(h+\eps)}<\mu(C_{ni}\cap D_j)<e^{-n(h-\eps)},
$$
and
$$
\tilde X_n=\{i\in X'_n\ |\ (i,j(i))\in Z_n\}.
$$
For $i\in\tilde X_n$
$$
\mu(C_{ni})\ge\mu(C_{ni}\cap D_{j(i)})>e^{-n(h+\eps)},
$$
and for $n$ large enough
$$
\mu(C_{ni})<{1\over1-\eps/2}\mu(C_{ni}\cap
D_{j(i)})<e^{n\eps}\mu(C_{ni}\cap D_{j(i)})<e^{-n(h-2\eps)}.
$$
Since for large $n$
$$
\sum_{(i,j)\in Z_n}\mu(C_{ni}\cap D_j)>1-\eps
$$
and
$$
\sum_{i\in X'_n}\mu(C_{ni}\cap D_{j(i)})>(1-\eps/2)\sum_{i\in
X'_n}\mu(C_{ni})>1-\eps,
$$
we have also
$$
\sum_{i\in\tilde X_n}\mu(C_{ni})\ge\sum_{i\in\tilde
X_n}\mu(C_{ni}\cap D_{j(i)})>1-2\eps.
$$
\endveriff

\noindent{\it Proof of Theorem~\ref{2.1}.} Let $C_n$ be a masa in
$N_n$. We may suppose that $C_n\subset C_{n+1}$. Then we have to
prove that the McMillan theorem holds for $\{C_n\}_n$, $H(\alpha)$
and any $\eps>0$ (note that if a projection $z_n'$ is chosen in
$C_n$ as in the statement of the McMillan theorem then we can
replace $1-z_n'$  by its central support $p_n$ in $N_n$ and let
$z_n = 1-p_n$).

Using Proposition~\ref{2.2}(i) choose $m\in\N$ such that
$$
\left|H(\alpha)-{1\over m}H(N_{m-p}\,,\alpha^m)\right|<\eps^2_1,
$$
where $\eps_1$ is as in Lemma~\ref{2.4}.  Let $A_m$ be the von
Neumann subalgebra of $\pi_\tau(A)''$ generated by
$\alpha^{jm}(N_{m-p})$, $j\in\Z$, $D$ a masa in $N_{m-p}$. Then
$\displaystyle D_n=\mathop{\bigvee}^{\left[{n\over
m}\right]-1}_{j=0}\alpha^{jm}(D)$ is a masa in
$\displaystyle\mathop{\bigvee}^{\left[{n\over
m}\right]-1}_{j=0}\alpha^{jm}(N_{m-p})$. If $\alpha^m$ was
ergodic, we could apply the classical McMillan theorem to
$\{D_n\}_n$ and then make use of Lemma~\ref{2.4} to conclude that
it holds also for $\{C_n\}_n$. Since $\alpha^m$ can be
non-ergodic, consider the fixed point algebra
$Z\subset\pi_\tau(A)''$ with respect to $\alpha^m$. Since $\alpha$
is asymptotically abelian and ergodic on $\pi_\tau(A)''$, $Z$ is a
finite dimensional subalgebra of the center of $\pi_\tau(A)''$.
Let $\{z_i\}_{i\in X}$ be the atoms of~$Z$. The automorphism
$\alpha$ acts transitively on the set of atoms, so the systems
$(Az_i,\tau_i,\alpha^m|_{Az_i})$ are pairwise conjugate, where
$\tau_i=|X|\tau|_{Az_i}$ and $|X|$ denotes the cardinality
of~$X$. Since ${1\over|X|}\sum_iH(\alpha^m|_{Az_i})=H(\alpha^m)$,
we conclude that $H(\alpha^m|_{Az_i})=mH(\alpha)$, so for
$H_i={1\over m} H(\alpha^m|_{A_mz_i})$ we have
$$
H_i\le{1\over m}H(\alpha^m|_{Az_i})=H(\alpha).
$$
On the other hand,
$$
{1\over|X|}\sum_iH_i={1\over m}H(\alpha^m|_{A_m\vee Z})={1\over
m}H(N_{m-p}\vee Z,\alpha^m)={1\over
m}H(N_{m-p}\,,\alpha^m)>H(\alpha)-\eps^2_1.
$$
Thus if $\tilde X=\{i\in X\ |\ H_i>H(\alpha)-\eps_1\}$ then
$$
H(\alpha)-\eps^2_1<{|\tilde
X|\over|X|}H(\alpha)+{|X\backslash\tilde
X|\over|X|}(H(\alpha)-\eps_1),
$$
whence
\begin{equation} \label{e2.2}
\sum_{i\notin\tilde X}\tau(z_i)={|X\backslash\tilde
X|\over|X|}<\eps_1<\eps.
\end{equation}
We have also
$$
H(\alpha)={1\over
m}H(\alpha^m|_{Az_i})\le\liminf_{n\to\infty}{H_{\tau_i}(C_nz_i)\over
n}
$$
and
$$
\lim_{n\to\infty}{1\over|X|}\sum_i{H_{\tau_i}(C_nz_i)\over
n}=\lim_{n\to\infty}{H(C_n\vee Z)\over
n}=\lim_{n\to\infty}{H(C_n)\over n}=H(\alpha),
$$
so
$$
\lim_{n\to\infty}{H_{\tau_i}(C_nz_i)\over n}=H(\alpha)\ \ \forall\
i\in X.
$$

For any $i\in X$ the automorphism $\alpha^m$ is ergodic on $Az_i$,
so by the classical result the McMillan theorem holds for
$\{D_nz_i\}_n$, $H_i$ and $\eps_1$. Then by Lemma~\ref{2.4}, for
$i\in\tilde X$, it holds for $\{C_nz_i\}_n$, $H(\alpha)$ and
$\eps$ (with respect to the trace $\tau_i$). By virtue of
(\ref{e2.2}) it holds also for $\{C_n\vee Z\}_n$, $H(\alpha)$ and
$2\eps$. Finally, by Lemma~\ref{2.5} the McMillan theorem holds
for $\{C_n\}_n$, $H(\alpha)$ and $4\eps$.
\endverif

The method used in the proof can be applied to prove the following
weak form of the McMillan theorem under more general assumptions.

\begin{theorem}
Let $(A,\tau,\alpha)$ be an asymptotically abelian system with
locality, and $\tau$ an ergodic trace. Then the entropy
$H(\alpha)$ of Connes and St{\o}rmer of the system coincides with
Voiculescu's completely positive approximation entropy
$hcpa_\tau(\alpha)$.
\end{theorem}

\noindent{\it Proof.} Let $N$ be a local subalgebra of $A$. It
suffices to prove that $hcpa_\tau(N,\alpha)\le H(\alpha)$. Then by
Proposition~\ref{2.2}(ii) it is enough to prove that
$hcpa_\tau(N_{m-p}\,,\alpha^m)\le mH(\alpha)$. Keep the notations
of the proof of Theorem~\ref{2.1}. By the classical McMillan
theorem we have $hcpa_{\tau_i}(N_{m-p}z_i\,,\alpha^m)
=H(N_{m-p}z_i\,,\alpha^m)$ (see the proof
of~\cite[Proposition~1.7]{V}). Since $H(N_{m-p}z_i\,,\alpha^m)\le
H(\alpha^m|_{Az_i})=mH(\alpha)$ and
$$
hcpa_\tau(N_{m-p}\,,\alpha^m)\le hcpa_\tau(N_{m-p}\vee
Z,\alpha^m)\le\max_i hcpa_{\tau_i}(N_{m-p}z_i\,,\alpha^m),
$$
we obtain the desired inequality.
\endverif

\begin{remarks}
\rm (i) The assumption that $N$ is a mean generator in
Theorem~\ref{2.1} is very close to being necessary. Indeed,
suppose $N$ satisfies the conclusion of Theorem. Then
$$
\log\rank N_n(1-z_n)<n(H(\alpha)+\eps).
$$
Thus
$$
H(N_n(1-z_n)\oplus\C z_n)\le\log\rank(N_n(1-z_n)\oplus\C z_n)
\le\log(e^{n(H(\alpha)+\eps)}+1),
$$
so
$$
\limsup_n{1\over n}H(N_n(1-z_n)\oplus\C z_n)\le
H(\alpha)+\eps\le\liminf_n{1\over n}H(N_n)+\eps.
$$
In particular, if the $N_n$'s don't grow too fast, e.g. if $\rank
N_n\le e^{Cn}$ for some $C>0$, then
$$
\limsup_n{1\over n}H(N_n)\le\limsup_n\left(C\tau(z_n)+{1\over n}
H(N_n(1-z_n)\oplus\C z_n)\right),
$$
and hence
$$
H(\alpha)=\lim_n{H(N_n)\over n}.
$$

\noindent(ii) It is not clear what the optimal assumptions are for
Theorem~\ref{2.1} to be true. The conclusion holds in several
cases when the dynamical system is not asymptotically abelian.
Such an example is that of a binary shift, see~\cite{PP}. Then we
are given a subset $X$ of $\N$ and the algebra $A(X)$ generated by
symmetries $\{s_n\}_{n\in\Z}$ satisfying commutation relations
$$
s_is_j=\cases{\ \ s_js_i,\ \ \hbox{if}\ \ |i-j|\notin X,\cr
-s_js_i,\ \ \hbox{if}\ \ |i-j|\in X.}
$$
The binary shift is the automorphism $\alpha$ defined by
$\alpha(s_n)=s_{n+1}$. Let $A_n=C^*(s_0,\ldots,s_{n-1})$. Then
$A_n\cong\Mat_{2^{d_n}}(\C)\otimes Z_n$, where $Z_n$ is the
diagonal in $\Mat_{2^{c_n}}(\C)$, $n=2d_n+c_n$, and thus any
minimal projection in $A_n$ has trace $2^{-d_n-c_n}$. We show
that the following three conditions are equivalent. By~\cite{GS1}
they are satisfied not only for asymptotically abelian systems.

(i) The conclusion of Theorem~\ref{2.1} holds for each algebra
$A_n$ (instead of $N$).

(ii) $H(\alpha)={1\over2}\log2$ and $\lim_n{c_n\over n}=0$.

(iii) $H(\alpha)=\lim_n{1\over n}H(A_n)$.

Indeed, the implication (i)$\Rightarrow$(iii) follows from the
previous remark.

Since $H(\alpha)\le\liminf_n{1\over n}H(A_n)={1\over2}\log2$, see
e.g. \cite{GS2}, (iii) implies $\lim_n{1\over
n}H(A_n)={1\over2}\log2$, hence by~\cite[Lemma 4.7]{GS2},
$\lim_n{c_n\over n}=0$. Thus (iii)$\Rightarrow$(ii).

If (ii) holds then ${1\over n}(d_n+c_n)\to{1\over 2}$, and so (ii)
implies (i).
\end{remarks}

\bigskip\bigskip

\section{Mean Generators} \label{3}

In this section we discuss several sufficient conditions for a
local algebra to be a mean generator.

\medskip

Our first result shows that it is enough to prove that at least
one algebra is a mean generator.

\begin{proposition} \label{3.1}
If $N$ is a mean generator for $A(N)$, $H(\alpha|_{A(N)})<\infty$
and $M$ is a subalgebra of $\vee^n_{-n}\alpha^i(N)$ for some $n$,
then $M$ is a mean generator for $A(M)$.
\end{proposition}

\noindent{\it Proof.} Without loss of generality we may suppose
that $M\subset N$. As above, let $p$ be such that $\alpha^j(N)$
commutes with $N$ for $|j|\ge p$. Let $D_n$ be a masa in
$M_{n-p}$, $C_n$ a masa in $N_{n-p}$ containing $D_n$. We have
$$
0\le H(D_n)-H(D_n,\alpha^n)=H(C_n)-(H(D_n,\alpha^n)+H(C_n|D_n))
\le H(C_n)-H(C_n,\alpha^n).
$$
Since ${1\over n}(H(C_n)-H(C_n,\alpha^n))\to0$ by assumption and
Proposition~\ref{2.2}, we conclude that
$$
{1\over n}(H(D_n)-H(D_n,\alpha^n))\to0,
$$
hence $M$ is a mean generator.
\endverif

We have used in the proof that if $H(\alpha|_{A(N)})<\infty$ then
$N$ is a mean generator if and only if
$$
\lim_{n\to\infty}{1\over n}(H(N_{n-p})-H(N_{n-p}\,,\alpha^n))=0.
$$
Using the following lemma we shall show that it suffices to check
this condition using any subalgebra of $N_n$ containing the
center $Z(N_n)$ instead of $N_n$.

\begin{lemma} \label{3.2}
Let $M_1,\ldots,M_n$ be commuting finite dimensional algebras,
$Z(M_i)\subset A_i\subset M_i$. Then
$$
H(\vee_iM_i)-H(\vee_iA_i)=\sum_i(H(M_i)-H(A_i)).
$$
\end{lemma}

\noindent{\it Proof.} It suffices to consider the case $n=2$. Let
$\{e_i\}_i$ be the atoms of a masa in $A_1$, $\{f_j\}_j$ the atoms
of a masa in $A_2$. Let $(M_1)_{e_i}$ be a factor of type
I$_{m_i}$, $(M_2)_{f_j}$ a factor of type I$_{n_j}$. Then
\begin{eqnarray*}
H(M_1\vee M_2)-H(A_1\vee A_2)
&=&-\sum_{i,j}\tau(e_if_j)\log{\tau(e_if_j)\over m_in_j}+
\sum_{i,j}\tau(e_if_j)\log\tau(e_if_j)\\
&=&\sum_i\tau(e_i)\log m_i+\sum_j\tau(f_j)\log n_j\\
&=&(H(M_1)-H(A_1))+(H(M_2)-H(A_2)).
\end{eqnarray*}
\endveriff

\begin{proposition} \label{3.3}
Let $N$ be a local algebra such that $H(\alpha|_{A(N)})<\infty$.
Let $\{M_n\}_n$ be a sequence of algebras such that
$Z(N_{n-p})\subset M_n\subset N_{n-p}$. Then $N$ is a mean
generator for $A(N)$ if and only if
$$
\lim_{n\to\infty}{1\over n}(H(M_n)-H(M_n,\alpha^n))=0.
$$
In particular, $N$ is a mean generator if $\lim_n{1\over
n}H(Z(N_n))=0$.
\end{proposition}

\noindent{\it Proof.} For any $m\in\N$, by Lemma~\ref{3.2}
$$
H(\vee^{m-1}_{j=0}\alpha^{jn}(N_{n-p}))-H(\vee^{m-1}_{j=0}\alpha^{jn}(M_n))
=m(H(N_{n-p})-H(M_n)),
$$
hence
$$
H(N_{n-p})-H(N_{n-p}\,,\alpha^n)=H(M_n)-H(M_n,\alpha^n),
$$
what gives the result.
\endverif

Finally, recall the following simple condition (see e.g.
\cite{Ch}).

\begin{proposition} \label{3.4}
Let $N$ be a local algebra. Suppose there exists $p\in\N$ such
that
$$\tau(ab)=\tau(a)\tau(b)
$$
for any $a\in\vee^0_{-\infty}\alpha^i(N)$ and
$b\in\vee^\infty_p\alpha^i(N)$. Then $N$ is a mean generator for
$A(N)$.
\end{proposition}

The last condition shows that Theorem~\ref{2.1} can be applied to
asymptotically abelian binary shifts and canonical shifts on
towers of relative commutants. On the other hand,
Proposition~\ref{3.1} allows to apply Theorem~\ref{2.1} to systems
arising from topological dynamics, which were considered
in~\cite[Section 5]{NS}. Indeed, if a local algebra has a masa
lying in the diagonal then it is a mean generator since all
computations are reduced to the abelian case. Then
Proposition~\ref{3.1} shows that any local algebra is a mean
generator.

\end{document}